

\magnification 1200
\hsize 13.2cm
\vsize 20cm
\parskip 3pt plus 1pt
\parindent 5mm

\def\\{\hfil\break}


\font\seventeenbf=cmbx10 at 17.28pt

\font\twelvebf=cmbx10 at 12pt
\font\eightbf=cmbx8
\font\sixbf=cmbx6

\font\eighti=cmmi8
\font\sixi=cmmi6

\font\eightrm=cmr8
\font\sixrm=cmr6

\font\eightsy=cmsy8
\font\sixsy=cmsy6

\font\eightit=cmti8
\font\eighttt=cmtt8
\font\eightsl=cmsl8

\font\seventeenbsy=cmbsy10 at 17.28pt

\font\twelvebsy=cmbsy10 at 12pt
\font\tenbsy=cmbsy10
\font\eightbsy=cmbsy8
\font\sevenbsy=cmbsy7
\font\sixbsy=cmbsy6
\font\fivebsy=cmbsy5

\font\tenmsa=msam10

\font\sevenmsa=msam7
\font\fivemsa=msam5
\newfam\msafam
  \textfont\msafam=\tenmsa
  \scriptfont\msafam=\sevenmsa
  \scriptscriptfont\msafam=\fivemsa

\font\tenmsb=msbm10
\font\eightmsb=msbm8
\font\sevenmsb=msbm7
\font\fivemsb=msbm5
\newfam\msbfam
  \textfont\msbfam=\tenmsb
  \scriptfont\msbfam=\sevenmsb
  \scriptscriptfont\msbfam=\fivemsb
\def\Bbb{\fam\msbfam\tenmsb}

\font\tenCal=eusm10
\font\sevenCal=eusm7
\font\fiveCal=eusm5
\newfam\Calfam
  \textfont\Calfam=\tenCal
  \scriptfont\Calfam=\sevenCal
  \scriptscriptfont\Calfam=\fiveCal

\font\teneuf=eusm10
\font\teneuf=eufm10
\font\seveneuf=eufm7
\font\fiveeuf=eufm5
\newfam\euffam
  \textfont\euffam=\teneuf
  \scriptfont\euffam=\seveneuf
  \scriptscriptfont\euffam=\fiveeuf

\font\seventeenbfit=cmmib10 at 17.28pt

\font\twelvebfit=cmmib10 at 12pt
\font\tenbfit=cmmib10
\font\eightbfit=cmmib8
\font\sevenbfit=cmmib7
\font\sixbfit=cmmib6
\font\fivebfit=cmmib5
\newfam\bfitfam
  \textfont\bfitfam=\tenbfit
  \scriptfont\bfitfam=\sevenbfit
  \scriptscriptfont\bfitfam=\fivebfit


\catcode`\@=11
\def\eightpoint{%
  \textfont0=\eightrm \scriptfont0=\sixrm \scriptscriptfont0=\fiverm
  \def\rm{\fam\z@\eightrm}%
  \textfont1=\eighti \scriptfont1=\sixi \scriptscriptfont1=\fivei
  \def\oldstyle{\fam\@ne\eighti}%
  \textfont2=\eightsy \scriptfont2=\sixsy \scriptscriptfont2=\fivesy
  \textfont\itfam=\eightit
  \def\it{\fam\itfam\eightit}%
  \textfont\slfam=\eightsl
  \def\sl{\fam\slfam\eightsl}%
  \textfont\bffam=\eightbf \scriptfont\bffam=\sixbf
  \scriptscriptfont\bffam=\fivebf
  \def\bf{\fam\bffam\eightbf}%
  \textfont\ttfam=\eighttt
  \def\tt{\fam\ttfam\eighttt}%
  \textfont\msbfam=\eightmsb
  \def\Bbb{\fam\msbfam\eightmsb}%
  \abovedisplayskip=9pt plus 2pt minus 6pt
  \abovedisplayshortskip=0pt plus 2pt
  \belowdisplayskip=9pt plus 2pt minus 6pt
  \belowdisplayshortskip=5pt plus 2pt minus 3pt
  \smallskipamount=2pt plus 1pt minus 1pt
  \medskipamount=4pt plus 2pt minus 1pt
  \bigskipamount=9pt plus 3pt minus 3pt
  \normalbaselineskip=9pt
  \setbox\strutbox=\hbox{\vrule height7pt depth2pt width0pt}%
  \let\bigf@ntpc=\eightrm \let\smallf@ntpc=\sixrm
  \normalbaselines\rm}
\catcode`\@=12

\def\eightpointbf{%
 \textfont0=\eightbf   \scriptfont0=\sixbf   \scriptscriptfont0=\fivebf
 \textfont1=\eightbfit \scriptfont1=\sixbfit \scriptscriptfont1=\fivebfit
 \textfont2=\eightbsy  \scriptfont2=\sixbsy  \scriptscriptfont2=\fivebsy
 \eightbf
 \baselineskip=10pt}

\def\tenpointbf{%
 \textfont0=\tenbf   \scriptfont0=\sevenbf   \scriptscriptfont0=\fivebf
 \textfont1=\tenbfit \scriptfont1=\sevenbfit \scriptscriptfont1=\fivebfit
 \textfont2=\tenbsy  \scriptfont2=\sevenbsy  \scriptscriptfont2=\fivebsy
 \tenbf}

\def\twelvepointbf{%
 \textfont0=\twelvebf   \scriptfont0=\eightbf   \scriptscriptfont0=\sixbf
 \textfont1=\twelvebfit \scriptfont1=\eightbfit \scriptscriptfont1=\sixbfit
 \textfont2=\twelvebsy  \scriptfont2=\eightbsy  \scriptscriptfont2=\sixbsy
 \twelvebf
 \baselineskip=14.4pt}

\def\seventeenpointbf{%
 \textfont0=\seventeenbf  \scriptfont0=\twelvebf  \scriptscriptfont0=\eightbf
 \textfont1=\seventeenbfit\scriptfont1=\twelvebfit\scriptscriptfont1=\eightbfit
 \textfont2=\seventeenbsy \scriptfont2=\twelvebsy \scriptscriptfont2=\eightbsy
 \seventeenbf
 \baselineskip=20.736pt}


\newdimen\srdim \srdim=\hsize
\newdimen\irdim \irdim=\hsize
\def\NOSECTREF#1{\noindent\hbox to \srdim{\null\dotfill ???(#1)}}
\def\SECTREF#1{\noindent\hbox to \srdim{\csname REF\romannumeral#1\endcsname}}
\def\INDREF#1{\noindent\hbox to \irdim{\csname IND\romannumeral#1\endcsname}}
\newlinechar=`\^^J
\def\openauxfile{
  \immediate\openin1\jobname.aux
  \ifeof1
  \message{^^JCAUTION\string: you MUST run TeX a second time^^J}
  \let\sectref=\NOSECTREF \let\indref=\NOSECTREF
  \else
  \input \jobname.aux
  \message{^^JCAUTION\string: if the file has just been modified you may
    have to run TeX twice^^J}
  \let\sectref=\SECTREF \let\indref=\INDREF
  \fi
  \message{to get correct page numbers displayed in Contents or Index
    Tables^^J}
  \immediate\openout1=\jobname.aux
  \let\END=\end \def\end{\immediate\closeout1\END}}

\newbox\titlebox   \setbox\titlebox\hbox{\hfil}
\newbox\sectionbox \setbox\sectionbox\hbox{\hfil}
\def\folio{\ifnum\pageno=1 \hfil \else \ifodd\pageno
           \hfil {\eightpoint\copy\sectionbox\kern8mm\number\pageno}\else
           {\eightpoint\number\pageno\kern8mm\copy\titlebox}\hfil \fi\fi}
\footline={\hfil}
\headline={\folio}

\def\titlerunning#1{\setbox\titlebox\hbox{\eightpoint #1}}
\def\title#1{\noindent\hfil$\smash{\hbox{\seventeenpointbf #1}}$\hfil
             \titlerunning{#1}\medskip}

\newcount\numbersection \numbersection=-1
\def\sectionrunning#1{\setbox\sectionbox\hbox{\eightpoint #1}
  \immediate\write1{\string\def \string\REF
      \romannumeral\numbersection \string{%
      \noexpand#1 \string\dotfill \space \number\pageno \string}}}
\def\section#1{%
  \par\vskip0.666cm\penalty -100
  \vbox{\baselineskip=14.4pt\noindent{{\twelvepointbf #1}}}
  \vskip2pt
  \penalty 500
  \advance\numbersection by 1
  \sectionrunning{#1}}

\def\subsection#1{%
  \par\vskip0.5cm\penalty -100
  \vbox{\noindent{{\tenpointbf #1}}}
  \vskip1pt
  \penalty 500}

\newcount\numberindex \numberindex=0
\def\index#1#2{%
  \advance\numberindex by 1
  \immediate\write1{\string\def \string\IND #1%
     \romannumeral\numberindex \string{%
     \noexpand#2 \string\dotfill \space \string\S \number\numbersection,
     p.\string\ \space\number\pageno \string}}}

\newdimen\itemindent \itemindent=\parindent

\def\item#1{\par\noindent\hangindent\itemindent%
            \rlap{#1}\kern\itemindent\ignorespaces}
\def\itemitem#1{\par\noindent\hangindent2\itemindent%
            \kern\itemindent\rlap{#1}\kern\itemindent\ignorespaces}
\def\itemitemitem#1{\par\noindent\hangindent3\itemindent%
            \kern2\itemindent\rlap{#1}\kern\itemindent\ignorespaces}

\long\def\claim#1|#2\endclaim{\par\vskip 5pt\noindent
{\tenpointbf #1.}\ {\it #2}\par\vskip 5pt}

\def\proof{\noindent{\it Proof}}

\def\today{\ifcase\month\or
January\or February\or March\or April\or May\or June\or July\or August\or
September\or October\or November\or December\fi \space\number\day,
\number\year}

\catcode`\@=11
\newcount\@tempcnta \newcount\@tempcntb
\def\timeofday{{%
\@tempcnta=\time \divide\@tempcnta by 60 \@tempcntb=\@tempcnta
\multiply\@tempcntb by -60 \advance\@tempcntb by \time
\ifnum\@tempcntb > 9 \number\@tempcnta:\number\@tempcntb
  \else\number\@tempcnta:0\number\@tempcntb\fi}}
\catcode`\@=12

\def\bibitem#1&#2&#3&#4&%
{\hangindent=0.8cm\hangafter=1
\noindent\rlap{\hbox{\eightpointbf #1}}\kern0.8cm{\rm #2}{\it #3}{\rm #4.}}


\def\bN{{\Bbb N}}

\def\bQ{{\Bbb Q}}
\def\bR{{\Bbb R}}
\def\bT{{\Bbb T}}



\def\square{{\hfill \hbox{
\vrule height 1.453ex  width 0.093ex  depth 0ex
\vrule height 1.5ex  width 1.3ex  depth -1.407ex\kern-0.1ex
\vrule height 1.453ex  width 0.093ex  depth 0ex\kern-1.35ex
\vrule height 0.093ex  width 1.3ex  depth 0ex}}}
\def\qed{\kern10pt$\square$}
\def\hexnbr#1{\ifnum#1<10 \number#1\else
 \ifnum#1=10 A\else\ifnum#1=11 B\else\ifnum#1=12 C\else
 \ifnum#1=13 D\else\ifnum#1=14 E\else\ifnum#1=15 F\fi\fi\fi\fi\fi\fi\fi}
\def\msatype{\hexnbr\msafam}
\def\msbtype{\hexnbr\msbfam}
\mathchardef\restriction="3\msatype16   
\mathchardef\boxsquare="3\msatype03
\mathchardef\preccurlyeq="3\msatype34
\mathchardef\compact="3\msatype62
\mathchardef\smallsetminus="2\msbtype72   
\mathchardef\subsetneq="3\msbtype28
\mathchardef\supsetneq="3\msbtype29
\mathchardef\leqslant="3\msatype36   
\mathchardef\geqslant="3\msatype3E   
\mathchardef\stimes="2\msatype02
\mathchardef\ltimes="2\msbtype6E
\mathchardef\rtimes="2\msbtype6F


\let\wt=\widetilde
\let\wh=\widehat
\let\text=\hbox
\def\buildo#1^#2{\mathop{#1}\limits^{#2}}
\def\buildu#1_#2{\mathop{#1}\limits_{#2}}
\def\ort{\mathop{\hbox{\kern1pt\vrule width4.0pt height0.4pt depth0pt
    \vrule width0.4pt height6.0pt depth0pt\kern3.5pt}}}

\def\vlra{\mathrel{\smash-}\joinrel\mathrel{\smash-}\joinrel%
    \kern-2pt\longrightarrow}
\def\srelbar{\vrule width0.6ex height0.65ex depth-0.55ex}
\def\merto{\mathrel{\srelbar\kern1.3pt\srelbar\kern1.3pt\srelbar
    \kern1.3pt\srelbar\kern-1ex\raise0.28ex\hbox{${\scriptscriptstyle>}$}}}


\def\Pic{\mathop{\rm Pic}\nolimits}



\long\def\InsertFig#1 #2 #3 #4\EndFig{\par
\hbox{\hskip #1mm$\vbox to#2mm{\vfil\special{"
(/home/demailly/psinputs/grlib.ps) run
#3}}#4$}}
\long\def\LabelTeX#1 #2 #3\ELTX{\rlap{\kern#1mm\raise#2mm\hbox{#3}}}


\itemindent = 7mm

\title{Numerical character of the effectivity} 
\title{of adjoint line bundles}

\titlerunning{}

\vskip10pt

\centerline {\tenrm Fr\'ed\'eric CAMPANA\quad Vincent KOZIARZ} 
\smallskip
\centerline {\tenrm Mihai P\u AUN}

\vskip20pt

\noindent{\bf Abstract. \it {In this note we show that, for any log-canonical pair $(X, \Delta)$, $K_X+ \Delta$ is $\bQ$-effective if its Chern class
contains an effective $\bQ$-divisor.
}}

%
\section{\S 0. Introduction} 

\noindent  For the notions of klt and log-canonical pairs, we refer to [13]. The main result we will prove is the following

\claim Theorem 1|Let $X$ be a smooth, connected (complex) projective manifold, and let $\Delta$ be an effective $\bQ$-divisor on $X$, such that the pair $(X, \Delta)$ is lc. Assume that there exists a line bundle $\rho$ on $X$ such that $c_1(\rho)= 0$, and such that
$H^0\big(X, m(K_X+ \Delta)+ \rho\big)\neq 0,$
for some integer $m$, divisible enough (i.e. such that all coefficients appearing in $\Delta$ become integral). 

Then $h^0\big(X, m^\prime(K_X+ \Delta)\big)\geq h^0\big(X, m(K_X+ \Delta)+ \rho\big)>0$, for some suitable multiple $m'$ of $m$.
\endclaim

In the first part of this note we treat the theorem for klt pairs, and in the second part we prove the result in full generality, by combining the techniques of proof in the klt case with the special case of a reduced boundary with simple normal crossings, established by Y. Kawamata in [12]. Finally, in the third part, we draw some immediate consequences of Theorem~1.

%

\noindent The special case $\Delta=0$ was settled in the article [5] by F. Campana and T. Peternell (with an appendix by M. Toma). The main steps of their approach are as follows: the case $m= 1$ is treated by using the fundamental result of 
C. Simpson in [17], and Serre duality. The general case was reduced to the $m=1$ 
case by means of ramified cover techniques, central in additivity results, due to H. Esnault and E. Viehweg. 

The main ingredient of the proof we present here is a generalization 
by N. Budur (cf. [4] and the references there) of C. Simpson's result, also based on the study of ramified covers \`a la Viehweg; this is discussed in Section 1.D.
We also introduce a simple and very effective technique for the reduction to the $m= 1$ case, see 1.C below. 

If the boundary $\Delta$ is reduced and has simple normal crossings, Theorem 1 has been recently established in [12]. The arguments invoked are the same as the ones in [5], with the notable exception that the mixed Hodge decomposition 
of P. Deligne for quasi-projective manifolds is used.

An easy consequence of Theorem 1 is that $\kappa(X,K_X+ \Delta)=0$, provided that the numerical dimension of $K_X+\Delta$ is equal to zero: this particular case of the abundance conjecture is well-known when $(X, \Delta)$ is klt, and is due to N.~Nakayama in [14]. However, the methods of [14] do not seem to apply in the more general context of the previous theorem. Notice that our result shows that the abundance conjecture is of numerical character. Also, O. Fujino 
suggested that the results proved here can be used in order to derive 
Corollary 3 (at the end of this text); we thank him for this.

When the boundary $\Delta$ has some ample part, Theorem 1
follows from a classical argument of V. Shokurov, which combines the Kawamata-Viehweg  vanishing theorem with the fact that the Euler characteristic of a line bundle is a topological invariant, see [16], [15].

Part of the results of [18] and [12] respectively by Y.-T. Siu and Y. Kawamata also rest on Simpson's result in an essential way. Notice that the general log-canonical version of Theorem 1 above is stated in the last section of [12] of Kawamata's text. Despite the apparences, the proof given there however seems to differ significantly from ours, given below in Section 2.

After posting the present text, we were informed by C. Hacon that a statement analogous to Lemma 1 below is formulated in [6], Theorem 3.2.

\bigskip

\section{\S 1. Proof of Theorem 1 in the klt case}
\noindent
First, we set $d:=h^0\big(X, m(K_X+ \Delta)+ \rho\big)>0$ and we denote by $F$ the fixed part of the linear system $\vert m(K_X+ \Delta)+ \rho\big\vert$ so that for any non-zero $u\in H^0\big(X, m(K_X+ \Delta)+ \rho\big)$, $[u=0]=F+M_u$ where $M_u$ is an effective divisor (the ``movable part'' of $[u=0]$). Once for all, we choose a generic $u_0$, meaning that $M_{u_0}$ is reduced and that the support of $M_{u_0}$ has no common component with the supports of $F$ and $\Delta$.

From now on, we denote by $(W_j)_{j\in J}$ the set of irreducible components appearing either in the support of $\displaystyle \Delta=\sum_{j\in J} \mu^j.W_j$, or of $[u_0=0]$. In particular, we can write $F= \sum_{j\in J} a^jW_j$ and $M_{u_0}=\sum_{j\in J} b^jW_j$.

\smallskip

We divide the proof of Theorem 1 in a few steps.

\medskip
\noindent{\bf 1.A} {\sl We can assume that the $(W_j)_{j\in J}$ are non-singular and have normal crossings.}

Indeed this is completely standard: we consider a log-resolution $\mu: \wh X\to X$ of $(X, \Delta)$
such that the proper transforms of $(W_j)_{j\in J}$, and the exceptional divisors are non-singular, and have normal crossings.
The change of variables formula reads as
$$\mu^*(K_X+ \Delta)+ E= K_{\wh X}+ \wh \Delta$$
where $E$ is effective and $\mu$-exceptional, $\wh \Delta$ is effective, and 
$(\wh X, \wh \Delta)$ is klt (see e.g. [8]).

Together with Hartogs principle, the above formula shows that the statement we want to prove is preserved by the modification $\mu$, i.e. it is sufficient to prove that some multiple of 
$K_{\wh X}+ \wh \Delta$ has at least $d$ linearly independent sections in order to conclude. Remark also that the assumptions concerning $u_0$ are not affected.
\hfill\qed

\noindent In what follows we will not change the notations, but we keep in mind that we have the transversality property 1.A.

\medskip
\noindent{\bf 1.B} {\sl We can assume that $\Delta$ and the zero divisor of $u_0$ have no common component}. 
\smallskip

Recall that we wrote  
$F= \sum_{j\in J} a^jW_j$. We have
$$m(\sum_{j\in J} \mu^jW_j)-\sum_{j\in J} a^jW_j=m(\sum_{j\in J_1} \mu^j_0W_j)-\sum_{j\in J_2} a^j_0W_j
$$
where $\displaystyle \mu^j_0:= \max\Big(\mu^j- {1\over m}a^j, 0\Big)$ and $\displaystyle a^j_0:= \max \big(a^j-m\mu^j, 0\big)$, hence the sets $J_1$ and $J_2$ corresponding to non-zero coefficients are disjoint.

Now, if $u\in H^0\big(X, m(K_X+ \Delta)+ \rho\big)$, we have

$$m(K_X+ \sum_{j\in J} \mu^jW_j)+\rho\sim \sum_{j\in J} a^jW_j+M_u$$
and thus 
$$m(K_X+ \sum_{j\in J_1} \mu^j_0W_j)+\rho\sim \sum_{j\in J_2} a^j_0W_j+M_u=:E_u\leqno (1)$$
where we denote by the symbol ``$\sim $" the linear equivalence of the bundles in question.
In particular, if $M_{u_0}=\sum_{j\in J_3} W_j$ (recall that $M_{u_0}$ is reduced) where $J_3$ is the set of non-zero coefficients, it follows from our choice of $u_0$ that the $J_i$ are pairwise disjoint, and the claim is proved. Indeed, if we are able to produce $d$ linearly independent sections of some multiple of the $\bQ$-bundle $K_X+ \sum_{j\in J_1} \mu^j_0W_j$, the conclusion follows.

 \hfill\qed

\smallskip

\noindent For the rest of this note, we replace the divisor $\Delta$ in Theorem 1 with $\sum_{j\in J_1} \mu^j_0W_j$
and we denote the divisor $E_{u_0}$ by $E$ so that

$$m(K_X+ \sum_{j\in J_1} \mu^j_0W_j)+\rho\sim \sum_{j\in J_2} a^j_0W_j+\sum_{j\in J_3}W_j=E.\leqno (1')$$
Observe that the supports of $E$ and $\Delta$ have no common component, and that their union is snc.

\medskip
 \noindent{\bf 1.C} {\sl Reduction to the case $m= 1$}
 
 \smallskip

\noindent We write the bundle $(1')$ in adjoint form, as follows:
$$m(K_X+\Delta)+\rho\sim K_X+\Delta+{m-1\over m}E +{1\over m}\rho.\leqno (2)$$
In order to simplify the writing, we introduce the next notations:
$$\Big[{m-1\over m}E\Big]:=\sum_{j\in J_2} \Big[{m-1\over m}a^j_0\Big]W_j,$$
and 
$$\Big\{{m-1\over m}E\Big\}:=\sum_{j\in J_2} \Big\{{m-1\over m}a^j_0\Big\}W_j+\sum_{j\in J_3} {m-1\over m}W_j;$$ 
here we denote by $\{x\}$ and $[x]$ respectively the fractional part and the integer part of the real number $x$. We therefore have the decomposition
$${m-1\over m}E= \Big[{m-1\over m}E\Big]+ \Big\{{m-1\over m}E\Big\}.$$

%
 
\smallskip

\noindent We subtract next the divisor $\big[{m-1\over m}E\big]$ from both sides of (2), and we define 
$$\Delta^+:=\Delta+{m-1\over m}E- \Big[{m-1\over m}E\Big]=\Delta+ \Big\{{m-1\over m}E\Big\}.$$
Henceforth, we have the identity
$$E- \Big[{m-1\over m}E\Big]\sim K_X+\Delta^++{1\over m}\rho.\leqno(3)$$
Moreover remark that by $(1)$, for any $u\in H^0\big(X, m(K_X+ \Delta)+ \rho\big)$ the divisor
$$E_u-\Big[{m-1\over m}E\Big]\sim K_X+\Delta^++{1\over m}\rho
$$
is also effective.

\smallskip
\noindent In conclusion, we can define a {\sl line bundle} $L_1$ associated to the  effective $\bQ$-divisor $\Delta^+$ such that:

\noindent $\bullet$ The pair $(X, \Delta^+)$ is klt;

\noindent $\bullet$ The adjoint bundle $\displaystyle K_X+ L_1+ {1\over m}\rho$ has $d$ linearly independent sections.


\medskip
\noindent{\bf 1.D} {\sl A lemma by N. Budur}
 
\smallskip

The following result could be extracted from N. Budur's article [4]; for the convenience of the reader, we include a direct proof as well. We denote by ${\rm Pic}^{\tau}(X)\subset {\rm Pic}(X)$ the subgroup of line bundles whose first Chern class is torsion.

\claim Lemma 1 ([4])|Let $X$ be a connected complex projective manifold, and $\Delta^+$ an effective $\bQ$-divisor on $X$, with simple normal crossings support, and $(X, \Delta^+)$ klt.  Assume also that $\Delta^+\sim L_1$, for some $L_1\in {\rm Pic}(X)$. 
For each integer $k\geq 0$, define $L_k:=kL_1-[k\Delta^+]\sim\{k\Delta^+\}$ (we remark that this is consistent with the previous assumption).

\noindent Then for each $k, i$ and $q$ the set 
$$V^{q}_i(f,L_k)=\{\lambda\in{\rm Pic}^\tau(X)\,:\,h^{q}(X,K_X+ L_k+ \lambda)\geq i\}
$$
is a finite union of torsion translates of complex subtori of ${\rm Pic}^0(X)$.
\endclaim

\proof. We write $\Delta^+=\sum_{j\in J}\alpha^jD_j$, and  let $N$ be the smallest positive integer such that $N\alpha^j\in\bN$ for all $j$. Then $NL_1$ has a section whose zero divisor is $\sum_{j\in J} N\alpha^jD_j$. We take the $N$-th root and normalize it in order to obtain a normal cyclic cover $\pi:\widetilde X\rightarrow X$ of order $N$. Let $\eta:\widehat X\rightarrow \widetilde X$ be a resolution and $f:=\pi\circ\eta$.
After fundamental results of Esnault-Viehweg [7], [8] we know that $\widetilde X$ has rational singularities hence $\eta_*{K}_{\widehat X}=K_{\widetilde X}$ and $R^i\eta_*{K}_{\widehat X}=0$ for all $i>0$ (see [13, Theorem 5.10] for example). Moreover,

$$\pi_*{K}_{\widetilde X}=K_X\otimes\bigoplus_{k=0}^{N-1}L_k,$$
and $R^i\pi_*{ K}_{\widetilde X}=0$ for all $i>0$, since $\pi$ is finite.

\smallskip

Finally, for any line bundle $\lambda$ on $X$ and any $q$, we have
$$H^{q}(\widehat X, K_{\widehat X}+ f^* \lambda)\simeq H^{q}(\widetilde X, K_{\widetilde X}+ \pi^* \lambda)\simeq\bigoplus_{k=0}^{N-1} H^{q}(X, K_X+ L_k+  \lambda)$$
by the Leray spectral sequence.

If we apply the result of Simpson [17], we know by Serre duality (although $\wh X$ might not be connected) that for any $i$
$$V^{q}_i(f)=\{\lambda\in{\rm Pic}^\tau(X)\,:\,h^{q}(\widehat X, K_{\widehat X}+ f^* \lambda)\geq i\}$$
is a finite union of torsion translates of complex subtori of ${\rm Pic}^0(X)$ (if $\wh X_1,\dots,\wh X_r$ are the connected components of $\wh X$, just write
$$V^{q}_i(f)=\bigcup_{i_1+\dots+i_r=i}\left[\bigcap_{\ell=1}^{r}V^q_{i_\ell}\Bigl(f_{|\wh X_{\ell}}\Bigr)\right]).$$
Now, an observation of Budur, Arapura, Simpson (cf. [4], see also [1], [2]) shows that each of the sets
$$V^{q}_i(f,L_k)=\{\lambda\in{\rm Pic}^\tau(X)\,:\,h^{q}(X, K_X+ L_k+ \lambda)\geq i\}
$$
has the same structure. Their argument goes as follows: let $V$ be an irreducible component of $V^{q}_i(f,L_k)\subset {\rm Pic}^\tau(X)$. For all $0\leq \ell<N$, let $\iota_\ell=\max\{p\,:\,V\subset V^{q}_p(f,L_\ell)\}$ and $I=\sum_{\ell=0}^{N-1} \iota_\ell$. Then $V$ is an irreducible analytic subset of $V^q_{I}(f)$ and it is an irreducible component of $\cap_{\ell=0}^{N-1}V^q_{\iota_\ell}(f,L_\ell)$ because $V^q_{\iota_k}(f,L_k)\subset V^q_i(f,L_k)$. Moreover,
$$V^q_{I}(f)=\bigcup_{\iota'_0+\dots+\iota'_{N-1}=I}\left[\bigcap_{\ell=0}^{N-1}V^q_{\iota'_\ell}(f,L_\ell)\right]
$$
and by construction, $V$ is not included in $\cap_{\ell=0}^{N-1}V^q_{\iota'_\ell}(f,L_\ell)$ if $\iota'\not=\iota$. Therefore, $V$ is an irreducible component of $V^q_{I}(f)$ for which Simpson's theorem applies, and the lemma is proved. \hfill\qed

\medskip

\noindent{\bf 1.E} {\sl End of the proof}

\smallskip

\noindent We follow next the original argument in [5]: by the second bullet at the end of Section~1.C 
 the bundle $\displaystyle K_X+\Delta^++{1\over m}\rho$ has $d=h^0(X,m(K_X+\Delta)+\rho)$ linearly independent sections, which means that
$${1\over m}\rho\in V^0_d(f,L_1)$$
(in the notations of Lemma 1).
By the results discussed in Section 1.D we infer the existence of
a torsion line bundle $\rho_{{\rm tor}}$ and of an element $T\in\bT$, where $\bT$  is a subtorus of $\Pic^0(X)$ such that
$$\rho= m(\rho_{{\rm tor}}+T).
$$
On the other hand, we also have $(1-m)T+ \rho_{{\rm tor}}\in V^0_d(f,L_1)$, 
which by definition implies that the bundle
$$K_X+ \Delta^++ (1-m)T+ \rho_{{\rm tor}}=K_X+\Delta^+-{m-1\over m}\rho+m\rho_{\rm tor}
$$
admits $d$ linearly independent sections. Now, the relation (2) shows that
$$K_X+\Delta^+-{m-1\over m}\rho\sim m(K_X+\Delta)-\Big[{m-1\over m}E\Big],$$
which implies that the line bundle
$$m(K_X+\Delta)-\Big[{m-1\over m}E\Big]+m\rho_{\rm tor}$$
has also $d$ linearly independent sections,
and Theorem 1 immediately follows.\hfill\qed

\medskip


%
%
%
\claim{Remark 1}|{\rm We observe that unlike in [5], we do not use a ramified cover of $X$ in order to reduce ourselves to the case $m= 1$. However, a ramified cover is used in 1.D in order to ``remove the boundary" $\Delta$. 
}

\endclaim

\medskip
\section{\S2 Proof of Theorem 1 in the log-canonical case} 
 
\medskip
\noindent In this section, our goal is to prove Theorem 1 in full generality.

The first remark is that the reductions performed in Sections 
 1.A-1.B still apply in the current lc setting; hence we can assume 
that:

\noindent $\bullet$ We have a decomposition $D=B+\Delta$, where the support of $D$ is snc, and $B=[D]$, so that $(X,\Delta)$ is klt.

\noindent $\bullet$ The union of the support of $D$ together with the support of the zero-locus $E$ of a chosen generic section $u_0\in H^0\big(X, m(K_X+D)+ \rho\big)$ is also a divisor which is snc.

\noindent $\bullet$ The divisors $E$ and $D$ have no common component.
 
 \medskip
\noindent We will show next that the proof of Theorem 1 is obtained as a consequence of 
two of its special cases: the case where $D=B$, treated by Y. Kawamata in [12], and (simple modifications of) the klt case treated in the previous section.

The arguments provided for the case where $D=B$ in [12] are parallel to the ones in [5], but in the mixed Hodge theoretic context of Deligne. We remark that although this result is stated in [12] only under the assumption that the numerical dimension of 
$K_X+D$ is equal to zero, the arguments given imply the general version. 
 
We shall reduce next to this case, by performing the constructions made in Section 1.C above (we simply replace $K_X$ with $K_X+ B$); we have
$$m(K_X+B+\Delta)+\rho\sim K_X+B+\Delta+{m-1\over m}E +{1\over m}\rho\leqno (4)$$
and
$$E-\Big[{m-1\over m}E\Big]\sim K_X+ B+ \Delta^+ +{1\over m}\rho\leqno(5)$$
with the same notations as above.

\medskip
Let $f= \pi\circ\eta:\widehat X\to X$ be the map associated to $N\Delta^+$ (cf. 1.D); recall that
$$f_*{K}_{\widehat X}=\bigoplus_{k=0}^{N-1}(K_X+ L_k),$$
and thus
$$f_*({K}_{\widehat X}+ f^*B+ {1\over m}f^*\rho)=\bigoplus_{k=0}^{N-1}(K_X+ B+ L_k+{1\over m}\rho).\leqno (6)$$
We will show in the lemma below that we have 
$$f_*({K}_{\widehat X}+f^*B+{1\over m}f^*\rho)=f_*({K}_{\widehat X}+(f^*B)_{\rm red}+{1\over m}f^*\rho).\leqno (7)$$
\noindent Granted the equality (7), the proof of Theorem 1 ends as follows.

By hypothesis, the bundle $K_X+B+L_1+{1\over m}\rho$ appearing in
the right-hand side of the equality (6) has $d$ linearly independent sections. We may assume that $(f^*B)_{\rm red}$ has snc support, from which we obtain by [12] that {\sl the set of $\rho$'s in ${\rm Pic}^{\tau}(X)$ for which  ${K}_{\widehat X}+(f^*B)_{\rm red}+f^*\rho$ has $d$ linearly independent sections is a torsion translate of a subtorus in ${\rm Pic}^{\tau}(X)$}. From the proof of  Lemma 1, we derive the same conclusion for each of the sets associated to the bundles $K_X+B+L_k$. The rest of the proof is the same as in the klt case. \hfill\qed

\medskip 

\noindent We prove now the equality (7).
\claim Lemma 2|With the above notations, we have
$$f_*\big({K}_{\widehat X}+f^*B+{1\over m}f^*\rho\big)=f_*\big({K}_{\widehat X}+(f^*B)_{\rm red}+{1\over m}f^*\rho\big).$$
\endclaim

\proof. It suffices to prove the lemma under the assumption that ${1\over m}\rho$ is trivial.  We first have an equality $\eta_* K_{\widehat X}= K_{\wt X}$, since the singularities of $\widehat X$ are rational.

We can now write $K_{\widehat X}+f^*B=K_{\widehat X }+\wh B+F$, with $\wh B=(f^*B)_{\rm red}$, and $F$ effective, contained  in the support of $f^*B$, and $\eta$-exceptional. By [13, Theorem 5.20], the pair $(\wt X,\pi^*B)$ is log-canonical. This means (by the very definition) that we have
$$K_{\widehat X}+\wh B\geq \eta^*\big(K_{\wt X}+ \pi^*B\big).$$
>From the preceding two observations, we infer that:
$$\eqalign{
 K_{\wt X}+ \pi^*B=&\eta_*(K_{\widehat X})+ \pi^*B= \eta_*(K_{\widehat X}+ f^*B)\geq \eta_*(K_{\widehat X}+\wh B)\geq \cr\geq & \eta_*(\eta^*(K_{\wt X}+ \pi^*B))=K_{\wt X}+ \pi^*B.\cr
}$$
The conclusion now follows by taking the direct images by $\pi$.\hfill\qed

\claim Remark 2|{\rm As we have already mentioned, Lemma 1 thus holds true (with the same proof) when $\Delta^+$ is, more generally, assumed to be only log-canonical.}
\endclaim

\medskip
\section{\S3 Some consequences of Theorem 1} 
 
\medskip

We start with an elementary remark.

\claim{Remark 3}|{\rm As a by-product of our proof, we obtain a very precise control of the zero set of the sections produced by Theorem 1. We refer e.g. to [15], Sections 1.G and 1.H for the relevance of this matter in the context of extension of twisted pluricanonical sections.

\noindent The set-up is as follows. We assume that there exist a set $J_0\subset J$ together with 
a set of rational numbers $0\leq d^j\leq \mu^j$ where $j\in J_0$, such that the effective divisor
$$D:= m\sum_{j\in J_0}d^jW_j\leq m\Delta$$
is contained in the fixed part $F$ of the linear system $|m(K_X+ \Delta)+ \rho|$.
Then we claim that the divisor 
$\displaystyle {m^\prime\over m}D$
is contained in the zero set of each of the sections in $H^0\big(X, m^\prime(K_X+ \Delta)\big)$ that we produce.
Indeed, the step 1.B of our proof consists in removing the common divisor between the boundary $\Delta$ and the fixed part $F$ (cf. the definition of 
$\mu^j_0$ and $a^j_0$). This divisor is bigger than $D$ so the claim trivially follows.
}

\endclaim
\medskip

\noindent As a consequence of Theorem 1, we obtain first the following statement
(which is analogous to results obtained in [5]).

\claim Corollary 1|Assume $(X,\Delta)$ is lc with rational coefficients. For every $\rho\in \Pic^{\tau}(X)$, we have
$\kappa(X,K_X+\Delta)\geq \kappa(X,K_X+\Delta+\rho)$.\endclaim

\proof. We have shown that the left member of the inequality is nonnegative
if so is the right member. Now assume that the right hand side is nonnegative, and consider the Moishezon-Iitaka fibrations
$f:X\to B$ and $g: X\to C$ of the $\bQ$-bundles $K_X+\Delta+\rho$ and
$K_X+\Delta$ respectively. These fibrations can be assumed to be regular
(by using  blow-ups of $X$ if necessary).

Let $Z$ be the general fibre of $g$; we denote by $\Delta_Z$ the restriction of 
$\Delta$ to $Z$. Then the pair $(Z,\Delta_Z)$ is still
lc, and $0=\kappa(Z,K_Z+\Delta_Z)\geq\kappa(Z,K_Z+\Delta_Z+\rho_Z)\geq 0$, as we see 
by using induction on the dimension of $X$ (if $X=Z$ the result is clear, by Theorem 1, because the conditions $\kappa\geq 0$ and $\kappa>0$ are characterised by an inequality on the number of linearly independent sections). But the preceding inequality shows that $f(Z)$ is zero-dimensional.

We conclude the existence of $h:C\to B$
such that $f=h\circ g$, and thus of the claimed inequality.
\hfill\qed
\medskip

\noindent We present next an $\bR$-version of Theorem 1; again, the motivation is  provided by the proof of the non-vanishing theorem in [3], [15]. 
 
Let $\displaystyle (W_j)_{j\in J_g\cup J_d}$ be a set of non-singular hypersurfaces of $X$ having normal crossings, where $J_g$ and $J_d$ are finite and disjoint sets. Denoting by ``$\equiv$'' the numerical equivalence of $\bR$-divisors, we assume that we have 
$$K_X+ \sum_{j\in J_g}\nu^jW_j\equiv \sum_{j\in J_d}\tau^jW_j\leqno (8)$$
where for each $j\in J_g$ we have $\nu^j\in ]0, 1]$, and for $j\in J_d$ we have 
$\tau^j> 0$. The numbers $(\nu, \tau)$ are not necessarily rationals. We state next the following direct consequence of Theorem 1.

\claim Corollary 2|Let $\eta> 0$ be a real number. Then for each $j\in J_g$ and $l\in J_d$ there exists 
a finite set of rational numbers $\displaystyle \Big({p^j_{k, \eta}\over q_\eta}, {r^l_{k, \eta}\over q_\eta}\Big)_{k=1,...,N_\eta}$ such that:
{
\itemindent 5.5mm
\item {\rm (a)} The vector $(\nu, \tau):= (\nu^j, \tau^l)$ is a convex combination of $\displaystyle (\nu_{k\eta}, \tau_{k\eta}):= \Big({p^j_{k, \eta}\over q_\eta}, {r^l_{k, \eta}\over q_\eta}\Big)$;

\item {\rm (b)} We have $|p^j_{k, \eta}- q_\eta\nu^j|\leq \eta$ and $|r^l_{k, \eta}- q_\eta\tau^l|\leq \eta$
for each $j, l, k, \eta$;
\smallskip

\item {\rm (c)} The bundle 
$$K_X+ \sum_{j\in J_g}{p^j_{k, \eta}\over q_\eta}W_j$$
is $\bQ$-effective. 

\noindent Hence, the bundle $K_X+ \sum_{j\in J_g}\nu^jW_j$ is $\bR$-linearly equivalent with an effective $\bR$-divisor.

}

\endclaim

\proof. We consider the set ${\cal A}\subset \bR^{|J_g|}\times \bR^{|J_d|}$ given by the couples $(x, y)$ such that 
$$K_X+ \sum_{j\in J_g}x^jW_j\equiv \sum_{j\in J_d}y^jW_j.\leqno (9)$$
We remark that ${\cal A}$ is non-empty, since by relation (8) it contains the point $(\nu, \tau)$. Also, it is an affine space defined over $\bQ$.

The upshot is that given $\eta> 0$, we can write $(\nu, \tau)$ as a convex combination of points $(\nu_{k,\eta}, \tau_{k, \eta})\in \cal A$ having 
rational coordinates, in such a way that the Dirichlet conditions in $\rm (b)$
are satisfied (see e.g. [15]).  

If $\eta\ll 1$, then as a consequence of $\rm (b)$ we infer that the coordinates of $\nu_{k,\eta}$
belong to $]0, 1]$, and that the coordinates of $\tau_{k, \eta}$ are positive. Hence the point 
$\rm (c)$ follows from Theorem 1. \hfill\qed
 
\medskip

\claim Remark 4|{\rm If the boundary divisor
$\Delta$ contains an ample part, then we can take $m^\prime:= m$ in Theorem 1 (see [15], Section 1.G). It would be particularly useful to have an analogous statement in our current setting.
}
\endclaim
 
 \noindent The following interesting corollary was brought to our attention by O. Fujino; it is a slightly more general version of a result due to S. Fukuda (see [10], in which $D$ is supposed to be semi-ample), appeared in discussions between O. Fujino and S. Fukuda.

\claim Corollary 3|Let $(X, \Delta)$ be a projective klt pair and let $D$ be a nef and 
abundant  $\Bbb Q$-Cartier $\Bbb Q$-divisor on $X$. 
Assume that $K_X+\Delta$ is numerically equivalent to $D$. 
Then $K_X+\Delta$ is semi-ample. 

\endclaim

 \proof. The {\sl nefness} of a line bundle is clearly a numerical property, hence 
 $K_X+\Delta$ is nef. 
 Our next claim is that we have the following sequence of relations
 $$\nu(K_X+\Delta)\geq \kappa(K_X+\Delta)\geq \kappa(D)=\nu(D)=\nu (K_X+\Delta).
 $$
 Indeed, the first inequality is valid for any nef bundle; the second one is the content of Corollary 1. The third relation above is a consequence of the fact that $D$ is nef and abundant, and the last one is due to the fact that the numerical dimension 
 $\nu$ of a nef line bundle only
 depends on its first Chern class.
 Thus, the corollary follows as a consequence of a result due to Kawamata in [11] (see also the version by O. Fujino in [9]). \hfill\qed

\medskip

\claim Remark 5|{\rm In the preceding corollary, if $(X, \Delta)$ is only lc, we cannot conclude (as Y. Gongyo pointed to us), since Y. Kawamata's theorem is no longer available. In dimension $4$, the conclusion nevertheless still holds, as shown in Y. Gongyo's article arXiv: 1005.2796.}
\endclaim


\vskip 15pt

\section{References}
\bigskip

{\eightpoint

\bibitem [1]&Arapura, D.:&\ Higgs line bundles, Green-Lazarsfeld sets, and maps of K\"ahler manifolds to
curves;& Bull.\ Amer.\ Math.\ Soc.\ 26 (1992), no. 2, 310-314&

\bibitem [2]&Arapura, D.:&\ Geometry of cohomology support loci for local systems. I;& J.\ Algebraic Geom.
6 (1997), no. 3, 563-597&

\bibitem [3]&Birkar, C., Cascini, P., Hacon, C., McKernan, J.:&\ Existence of minimal models for varieties of log general type;&\ J. Amer. Math. Soc. 23 (2010), 405-468&

\bibitem [4]&Budur, N.:&\ Unitary local systems, multiplier ideals, and polynomial periodicity of Hodge numbers;& arXiv:math/0610382, to appear in Adv.\ Math.& 

\bibitem [5]&Campana, F., Peternell, T., Toma, M.:&\ Geometric stability of the cotangent bundle and the universal cover of a projective manifold;& arXiv:math/0405093, to appear in Bull.\ Soc.\ Math.\ France& 

\bibitem [6]&Chen, J., Hacon, C.:&\ On the irregularity of the image of the Iitaka fibration;&\ Comm. in Algebra Vol. 32, No. 1, pp. 203-215 (2004)&

\bibitem [7]&Esnault, H., Viehweg, E.:&\ Logarithmic de Rham complexes and vanishing theorems;& Invent.\ Math.\ 86 (1986), no. 1, 161-194&

\bibitem [8]&Esnault, H., Viehweg, E.:&\ Lectures on vanishing theorems;& DMV Seminar, 20.
Birkh\"auser Verlag, Basel, 1992&

\bibitem [9] &Fujino, O.:& On Kawamata's theorem; &arXiv:0910.1156&

\bibitem [10] &Fukuda, S.:& An elementary semi-ampleness result for log-canonical divisors;& arXiv:1003,1388&

\bibitem [11]&Kawamata, Y.:& Pluricanonical systems on minimal algebraic varieties ;& Inv. Math. 79 (1985), 567-588&

\bibitem [12]&Kawamata, Y.:&\ On the abundance theorem in the case $\nu=0$;& arXiv:1002.2682&

\bibitem[13] &Koll\`ar, J., Mori, S. : & Birational geometry of Algebraic Varieties. Cambridge university tracts 134.& Cambridge university press (1998)&

\bibitem [14]&Nakayama, N.:&\ Zariski decomposition and abundance;&  MSJ Memoirs  {\bf 14}, Tokyo (2004)&

\bibitem [15]&P\u aun, M.:&\ Relative critical exponents, non-vanishing and metrics with minimal singularities;& arXiv:0807.3109&

\bibitem [16]&Shokurov, V.:&\ A non-vanishing theorem;&\ Izv. Akad. Nauk SSSR
(49) 1985&

\bibitem [17]&Simpson, C.:&\ Subspaces of moduli spaces of rank one local systems;& Ann.\ Sci.\ ENS. (4) 26 (1993), no. 3, 361-401&
 
\bibitem [18]&Siu, Y.-T.:&\ Abundance conjecture;& arXiv:0912.0576&

}

\bigskip
\noindent
\bigskip\bigskip
{\parindent=0cm
F. Campana, V. Koziarz, M. P\u aun\\
Institut Elie Cartan, Universit\'e
Henri Poincar\'e, B. P. 70239, F-54506 Vand\oe uvre-l\`es-Nancy Cedex, France\\
E-mails: campana@iecn.u-nancy.fr, koziarz@iecn.u-nancy.fr, paun@iecn.u-nancy.fr

}

\end